\documentclass[15pt,a4paper]{article}
\usepackage{amsmath,amssymb,amsthm}

\numberwithin{equation}{section}

\theoremstyle{definition}

\newtheorem{Rem}{Remark}[section]

\usepackage{amsmath, amsthm, amsfonts, amssymb}
\usepackage{mathptmx}
\usepackage{graphicx}
\usepackage{fancybox}
\usepackage{color}

\Large
\def\supp{\operatorname{supp}}

\def\div{\operatorname{div}}
\def\a{\alpha}
\def\e{\varepsilon}

\newtheorem{thm}{Theorem}[section]

\newtheorem{df}{Definition}[section]

\newtheorem{lem}{Lemma}[section]
\newtheorem{cor}{Corollary}[section]
\newtheorem{prop}{Proposition}[section]
\newtheorem{rem}{Remark}[section]
\numberwithin{equation}{section}

\newcommand{\wt}[1]{\widetilde{#1}}

\newcommand{\wbs}[1]{\,\overline{\!{#1}}}

\newcommand{\wbml}[1]{\;\!\overline{\:\!\!{#1}\;\!\!}}

\newcommand{\ch}[1]{\chi\kern-.05em\lower1ex\hbox{$\scriptstyle{#1}$}}
\newcommand{\chb}[1]
{\wbml{\chi}\kern-.02em\lower1ex\hbox{$\scriptstyle{#1}$}\;\!}

\newcommand{\sgn}{{\rm{sgn}}\,}

\renewcommand{\div}{{\rm{div}}\,}

\newcommand{\Rn}{{\bf R}\kern-0.08em\lower-0.75ex\hbox{$\scriptscriptstyle n$}}
\newcommand{\Rd}[1]{{\bf R}\kern-0.1em\lower-0.75ex\hbox{$\scriptscriptstyle{#1}$}}
\newcommand{\Sd}{S\kern0.0em\lower-0.75ex\hbox{$\scriptscriptstyle d\;\!\!-\:\!\!1$}}
\newcommand{\Sn}{S\kern0.0em\lower-0.75ex\hbox{$\scriptscriptstyle n\;\!\!-\:\!\!1$}}
\newcommand{\Cn}{{\bf C}\kern-0.01em\lower-0.75ex\hbox{$\scriptscriptstyle n$}}
\newcommand{\el}[1]{{\ell}\kern0.02em\lower-0.75ex\hbox{$\scriptscriptstyle{#1}$}}
\newcommand{\EL}[1]{{L}\kern0.0em\lower-0.75ex\hbox{$\scriptscriptstyle{#1}$}}
\newcommand{\EA}[1]{{A}\kern0.0em\lower-0.75ex\hbox{$\scriptscriptstyle{#1}$}}
\newcommand{\Tst}{T\kern0.0em\lower-0.58ex\hbox{$\scriptscriptstyle *$}}

\newcommand{\wtM}[1]{{\:\wt{\;\!\!\!M\;\!\!}}
\kern0.1em\lower-0.9ex\hbox{$\scriptstyle{#1}$}}
\newcommand{\wbM}[1]{{\;\!\wbs{\:\!\!M\:\!\!}}
\kern0.14em\lower-0.9ex\hbox{$\scriptstyle{#1}$}}

\title{Kato's inequalities up to  the  boundary  for a quasilinear elliptic operator}

\author{Toshio Horiuchi, Peter Kumlin }
\date{}

\begin{document}
\maketitle

\begin{abstract}
Let  $\Omega$ be a bounded  smooth domain of $ \mathbb{R}^{N}$. By $\Delta_p $ with $1<p<\infty$ we denote $p$-Laplacian.
We prove  that  if $\Delta_p u$ is a finite measure in $\Omega$, then under suitable assumptions on $u$, $\Delta_p u^+ $ is also a finite measure  in $\Omega$ up to the boundary $\partial \Omega$. 
\footnote{ 
This research was partially supported by Grant-in-Aid for Scientific Research (No. 16K05189) and (No. 15H03621).}
\end{abstract}

\section{Introduction}
Let $\Omega$ be a bounded smooth domain of $ \mathbb{R}^{N}$.  By $\Delta_p $ for $p\in (1,+\infty)$ we denote $p$-Laplacian.
The classical   Kato's inequality for  a Laplacian in \cite {d} asserts that  given any function $u\in L^1_{\rm loc} (\Omega)$ such  that  $\Delta u\in L^1_{\rm loc} (\Omega)$, then 
$\Delta (u^+) $  is a Radon measure and  the  following holds:
\begin{equation} \Delta (u^+) \ge \chi_{[u\ge 0]} \Delta u \qquad \mbox{  in  } D'(\Omega), \label{Kato}
\end{equation}
where $u^+= \max[u,0]$. 
In \cite{BP1, BP2},
    H.Brezis and A.Ponce intensively  studied Kato's inequalities  with $\Delta u$ being a Radon measure and established 
  the strong maximum principle, the improved Kato's inequality and the inverse maximum principle (See also \cite{BMP, DP}).
  Then, 
 in \cite{LHkato,LH} Kato's inequality was further  studied for $\Delta_p u$ with $p\in (1,\infty)$ and 
 most of the  counter-parts were established under the  assumption  that $u$ is  admissible in $W^{1,p^*}_{\rm loc}(\Omega)$, where  $p^* = \max( 1, p-1)$.  
For the  admissibility in  $W^{1,p^*}_{\rm loc}(\Omega)$, see Definition \ref{df3} in Appendix and see also  \cite{LHkato2}.  
  We  remark  that when $p=2$, the  notion of  admissibility becomes  trivial.
On the other hand,  H.Brezis and  A. Ponce  in \cite{BP3} and A. Ancona  in \cite{A} studied Kato's inequality (\ref{Kato}) up to the boundary for  $p=2$.@\par
  The purpose in  the present paper  is to study  Kato's inequality for  $\Delta_p  $     up to the boundary of  $\Omega$. As a result, we will show  that $\Delta_p u^+$  is also a finite measure under suitable  assumptions on $u$.
  In these arguments  it is crucial to introduce  a class ${\Bbb X}_p$ in Definition {\ref{1.1}},  which was originally introduced in Brezis, Ponce \cite{BP3} for $\Delta$,
   and to  use  effectively a  notion of    admissibility  in ${\Bbb X}_p$  for  $\Delta _p$.
   \begin{df} \label{1.1}
 We say $u\in {\Bbb X}_p$ if 
 $u\in W^{1,p^*}(\Omega)$  and if there exists  a constant $C>0$ such that 
 \begin{equation}\left|  \int_\Omega |\nabla u|^{p-2} \nabla  u\cdot \nabla \varphi \right| \le C ||\varphi||_{L^\infty(\Omega)}, \qquad \mbox{ for  any }  \varphi \in C^1(\overline\Omega),
 \end{equation}
 in which case we set 
 \begin{equation} [ u]_{{\Bbb X}_p}= \sup_{\substack{\psi \in C^{1}(\bar{\Omega}) \\ \|{\psi}\|_{L^{\infty}} \le 1}} 
  \int_{\Omega}|\nabla u |^{p-2}\nabla{u}\cdot \nabla{\psi}.
 \end{equation} 
\end{df}
If $u \in{{\Bbb X}_p}$, then there exists a unique  bounded linear functional $T \in [C(\bar{\Omega})]^{*} = \mathcal{M}_b(\bar{\Omega})$
such  that 
\begin{equation*}
\langle{T,\psi}\rangle =  \int_{\Omega}|\nabla u |^{p-2}\nabla{u}\cdot \nabla{\psi} \quad (\forall{\psi} \in C^{1}(\bar{\Omega})).
\end{equation*}
On the other hand, by the Riesz Representation Theorem any $T\in \mathcal{M}_b(\bar{\Omega})$ admits 
a unique decomposition 
\begin{equation*}
\langle{T,\psi}\rangle = \int_{\partial\Omega} \psi~d\nu + \int_{\Omega} \psi~d\mu \quad (\forall{\psi} \in C(\bar{\Omega})),
\end{equation*}
where 
$\mu \in \mathcal{M}_b(\Omega)$ and $\nu \in \mathcal{M}_b(\partial\Omega)$.
By  
$\mathcal{M}_b(\Omega)$   and $\mathcal{M}_b(\partial\Omega)$  we denote 
the space  of all bounded measures in $\Omega$  and $\partial\Omega$,  equipped with the standard  norms
$\|\cdot \|_{\mathcal{M}_b(\Omega)}$  and $\|\cdot \|_{\mathcal{M}_b(\partial\Omega)}$ respectively. 
 We remark  that  measures in $\mathcal{M}_b(\Omega)$ are identified  with measures in ${\Omega}$ which do not charge $\partial\Omega$.
More precisely we  have
\begin{align*}
 || \mu||_{\mathcal{M}_b(\Omega)} &= \sup \left\{ \int_\Omega \varphi\,d\mu;  \varphi\in C_0(\bar\Omega) \mbox{  and   }
 ||\varphi||_{L^\infty(\Omega)}\le 1\right\}, 
\end{align*}
where  by $C_0(\bar\Omega)$ we denote  the space of all continuous functions on $\bar\Omega $ vanishing on $\partial\Omega$.
On  the other hand
$\mathcal{M}(\Omega)$  denotes  the space  of all Radon  measures in $\Omega$. In other words $\mu \in \mathcal{M}(\Omega)$ if and only if, for every $\omega\subset\subset \Omega$, there is $C_\omega>0$ such  that  $| \int_\Omega \varphi\,d\mu| \le C_\omega ||\varphi||_\infty$ for all $\varphi \in C_0 (\overline\omega)$.
When $u \in{{\Bbb X}_p}$, we will denote 
\begin{equation*}
\mu = -\Delta_p u ~ , \quad \nu =|\nabla u|^{p-2} \frac{\partial{u}}{\partial{n}},
\end{equation*} where $n$  denotes the outer normal.
In this paper, 
for $u \in{{\Bbb X}_p}$  we  always use  the notations $\Delta_p u$  and  $|\nabla u|^{p-2} \frac{\partial{u}}{\partial{n}}$ in the above  sense.
Hence  if 
$u \in{{\Bbb X}_p}$, then  we have 
\begin{equation*}
\int_{\Omega} |\nabla u|^{p-2}\nabla{u}\cdot \nabla{\psi} = \int_{\partial\Omega} \psi |\nabla u|^{p-2}\frac{\partial{u}}{\partial{n}} - \int_{\Omega} \psi \Delta_p u \quad (\forall{\psi} \in C^{1}(\bar{\Omega})).
\end{equation*}
It follows  from Theorem \ref{thm0} that for  every  
$u\in {{\Bbb X}_p}$\begin{equation*}
[u]_{{{\Bbb X}_p}} =  \int_{\Omega} |{\Delta_p u}| + \int_{\partial\Omega}|\nabla u|^{p-2} \Bigl{|}{\frac{\partial{u}}{\partial{n}}}\Bigr{|}
\end{equation*}
and if
$u$ is  admissible  in $\Bbb X_p$, then    $[u]_{{{\Bbb X}_p}} = 0$ if  and only if 
$u = const. $ in $\Omega$. 


%
%

\section{Preliminaries:
 Admissibilities in  ${\Bbb X}_p$ and $W^{1,p^*}_{0}(\Omega)$}
  We will work with  the  standard Sobolev spaces; 
  $ W^{1,p}(\Omega)$ and $W^{1,p}_0(\Omega)$,  where
 the  space $W^{1,p}(\Omega)$  is  equipped with  the norm
\begin{equation} || u||_{W^{1,p}(\Omega)}=|| |\nabla u| ||_{ L^p(\Omega)}+ ||u||_{ L^p(\Omega)},
\end{equation} and  by
 $W^{1,p}_0(\Omega)$ we denote the  completion of $C_c^\infty(\Omega)$ in the norm $|| \cdot ||_{W^{1,p}(\Omega)}$.
Now  we  introduce   two admissiblities  for $\Delta_p$
 to deal with Kato's  inequalities up to  the boundary.
  We note  that  theses notions become  trivial  if  $p=2$ and 
  a local version  was already introduced in  \cite{LH}.
  \begin{df}\label{df1.2} {\bf  $($Admissibility    in ${\Bbb X}_p$ $)$}
  Let  $ 1<p<\infty $ and $p^* =\max[1,p-1]$.
 A function $u$ is said to be  admissible  in ${\Bbb X}_p$  if   $u \in {\Bbb X}_p$ and   
 there exists a sequence $\{ u_k\}_{k=1}^\infty \subset W^{1,p}(\Omega) \cap L^\infty(\Omega)$  such  that: 
 \begin{enumerate}
    \item  $u_k\to u $ a.e. in $\Omega$ and $\,  u_k\to u $ in $ W^{1,p^*}(\Omega)$ as $k \to \infty$.
 \item $\Delta_p u_k\in L^1(\Omega)\,$  and   $|\nabla u_k|^{p-2}\frac{\partial u_k}{\partial n} \in L^1(\partial \Omega)$ 
 $(k=1,2,\cdots)$ and 
 \begin{align} & \sup_{k}||\Delta_p u_k||_{\mathcal M_b(\Omega)}= \sup_{k}\int_\Omega |\Delta_p u_k|<\infty
  \label{ub} \\
 & \sup_{k}\left| \left| |\nabla u_k|^{p-2}\frac{\partial u_k}{\partial n} \right| \right|_{\mathcal M_b(\partial\Omega)}=
  \sup_{k} \int_{\partial\Omega} \left|  |\nabla u_k|^{p-2}\frac{\partial u_k}{\partial n} \right| 
   <\infty. \label{ub2}
 \end{align}
  \end{enumerate}
  \end{df}
         \begin{df}\label{df1.1} {\bf  $($Admissibility   in $W^{1,p^*}_{0}(\Omega)$  $) $}  Let  $ 1<p<\infty $ and $p^* =\max[1,p-1]$.
 A function $u$ is said to be  admissible  in $W^{1,p^*}_{0}(\Omega)$   if   $u \in W^{1,p^*}_{0}(\Omega)$,  $\Delta_p u \in \mathcal{M}_b(\Omega)$ and   
 there exists a sequence $\{ u_k\}_{n=1}^\infty \subset W^{1,p}_{0}(\Omega) \cap L^\infty(\Omega)$  such  that: 
 \begin{enumerate}
    \item  $u_k\to u $ a.e. in $\Omega$ and $\,  u_k\to u $ in $ W^{1,p^*}_{0}(\Omega)$ as $k \to \infty$.
 \item $\Delta_p u_k\in L^1(\Omega)\,$  $(k=1,2,\cdots)$ and 
 \begin{equation} \sup_{k}||\Delta_p u_k||_{\mathcal M_b(\Omega)}= \sup_{k}\int_\Omega |\Delta_p u_k|<\infty. \label{ub}
 \end{equation}
  \end{enumerate}
  \end{df}
  Roughly speaking, if $u$ is  admissible in one of  these definitions, then  $u$ can be approximated by a sequence of good functions not only in the sense of the distributions but also in the sense of  measures.   
  Moreover it  is  possible to approximate    $u$ by  a sequence  of $C^1$-functions provided  that $u$  is  admissible. In fact 
in  Proposition \ref{prop1.1}  in  Appendix we collect such  nice properties of admissible functions  together with a local version of  the admissibility in  $W^{1,p^*}_{\rm loc}(\Omega)$. In the subsequent we  describe more remarks.
 \begin{rem}\label{remark2.1}
 \begin{enumerate}
  \item For  a general 
  class of uniformly  elliptic  operators with a divergence  form,  one can  define  the admissibility and establish  similar results  in parallel to  the present paper (c.f. \cite{LHkato2}).
  Further it  is  possible  to construct non-admissible functions in such  cases. 
When  $p=2$,  the existence of pathological    solution, which  is  non-admissible,  was initially shown by J Serrin in  the famous paper   \cite{Sj} (See also  \cite{JMS}).
 \item
  If  $u \in  W^{1,p^*}_{\rm loc}(\Omega) $,  then $\Delta_p u$, $\Delta_p( u^+ )$  and $\Delta_p( u^-) $  are well-defined in  $D'(\Omega)$. Let $\{u_k\}$ be the  sequence  in one of the  definitions.
It  follows  from the condition 1 that 
 $\Delta_p u_k =\Delta_p (u_k^+ )- \Delta_p( u_k^-)$  and  $\Delta_p u_k\to \Delta_p u $ ( i.e. $\, \Delta_p (u_k^{\pm}) \to \Delta_p( u^{\pm} )$ )   in  $D'(\Omega)$ as  $k\to\infty$. 
  Moreover, it  follows from the  condition 2 and  the  weak compactness of  measures that   we   have $
\Delta_p u_k \to \Delta_p u  $ ( i.e.  $  \Delta_p (u_k^{\pm})  \to \Delta_p (u^{\pm} )$ )   in  the   sense of  measures  as  $n\to\infty$.
(In the case of Definition \ref{df1.2},  $|\nabla u_k|^{p-2}\frac{\partial u_k}{\partial n}\to |\nabla u|^{p-2}\frac{\partial u}{\partial n}$ in  the   sense of  measures   as  well.)
Therefore   if  $u$ is admissible,  then $u^+$  and $u^-$  are so  as well. 

  \item
Let  $\Omega$  be  a unit ball $B_1(0)$  of {\bf $R^{\it N}$}.
 Let  $u=|x|^{\a} -1$ for $ \a = (p-N)/(p-1)$ and $p\in (1,N)$. 
 Then $u$ satisfies 
$$\Delta_p u= \a|\a|^{p-2} c_N\delta \qquad \text{ in } D'(\Omega),$$
where $\delta $ denotes a Dirac mass and $c_N$ denotes the  surface area of  the $N$-dimensional  unit ball $B_1$. Then $u$ is  admissible  in $ W^{1,p^\ast}_0(\Omega) $ if  $p\in (2-1/N,N)$ with $N\ge2$.  
 We  note  that when $1<p<2-\frac 1N $, $u$ is not admissible but  regarded as a renormalized solution. For  the detail see \cite{BBGGPV,BGO,Maso, Maeda, MMOP}
   \end{enumerate}
   \end{rem}

\section{Main results}
Given $M > 0$, we denote a truncation function $T_M$: {\boldmath $R$} $\to$ {\boldmath $R$} by
\begin{align}
T_M(s)=\max[-M, \min[M,s]]. \label{tran}
\end{align}

\begin{thm} \label{thm0}  If  $u\in {{\Bbb X}_p}$,  then we have:

\begin{enumerate}
\item

\begin{equation}
[u]_{{{\Bbb X}_p}} =  \int_{\Omega} |{\Delta_p u}| + \int_{\partial\Omega}|\nabla u|^{p-2} \Bigl{|}{\frac{\partial{u}}{\partial{n}}}\Bigr{|}. \label{seminorm}
\end{equation}
\item If $u$ is  admissible  in $\Bbb X_p$,  then  for every  $M>0$ ~$T_M u\in W^{1,p}(\Omega)$   and  we  have 
\begin{equation}  \int_\Omega\left| \nabla T_M(u)\right|^p \le M [u]_{{{\Bbb X}_p}}.
\end{equation}
\item
If $u$ is  admissible  in $\Bbb X_p$, then    $[u]_{{{\Bbb X}_p}} = 0$ if  and only if 
$u = const. ~ in ~ \Omega$.
%
%

\end{enumerate}

\end{thm}

\begin{thm}\label{thm1}
If $u $  is  admissible in ${{\Bbb X}_p}$,  then $u^{+} \in{{\Bbb X}_p}$ and we have
\begin{equation} 
[u^{+}]_{{{\Bbb X}_p}} \le [u]_{{{\Bbb X}_p}}.
\end{equation}
\end{thm}

\begin{thm} \label{thm2}
Assume that $u $  is admissible in $ W^{1,p^*}_{0}(\Omega)$.  Then we  have the  followings:  
\begin{enumerate}
\item   $u $  is  admissible in ${{\Bbb X}_p}$~( hence A$u^{+} \in{{\Bbb X}_p}$ ).
\item  \begin{equation} 
\int_{\Omega} |\Delta_p u^{+}| \le \int_\Omega |\Delta_p u|. \label{3.5}
\end{equation}

\end{enumerate}

\Rem If 
$u$  does not vanish on $\partial{\Omega}$, then  the  assertion (\ref{3.5}) fails. To  see this it suffices to take a linear function $u$.
\end{thm}

\begin{thm} \label{thm3} Assume that 
$u$ is  admissible in ${{\Bbb X}_p}$. Moreover  assume that $\Delta_p u \in L^{1}(\Omega),~ |\nabla u|^{p-2}\frac{\partial{u}}{\partial{n}} \in L^{1}(\partial{\Omega})$.  Then 
\begin{equation} 
\int_\Omega  |\nabla u|^{p-2}\nabla{u^{+}}\cdot \nabla{\psi} \le \int_{\partial{\Omega}} H \psi - \int_{\Omega} G \psi \quad (\forall{\psi} \in C^{1}(\bar{\Omega}) , \psi \ge 0 ~ in ~ \Omega)
\end{equation}
Here $G \in L^{1}(\Omega)$  and $H \in L^{1}(\partial{\Omega})$ are  given by
\begin{align} 
G = \begin{cases}
		\Delta_p u &on ~ [u > 0] \\
		0 &on ~ [u \le 0]
	\end{cases} ~ , \quad
H = \begin{cases}
		 |\nabla u|^{p-2}\frac{\partial{u}}{\partial{n}} &on ~ [u > 0] \\
		0 &on ~ [u < 0] \\
		\min \{  |\nabla u|^{p-2}\frac{\partial{u}}{\partial{n}},0 \} &on ~ [u = 0]
	\end{cases}
\end{align}
Thus, we have 
\begin{align} 
\begin{cases}
		\Delta_p u^{+} \ge G &in ~ \Omega \\
		~ |\nabla u|^{p-2}\frac{\partial{u^{+}}}{\partial{n}} \le H &on ~ \partial{\Omega}.
\end{cases}
\end{align}

\end{thm}

\subsection{Proof of Theorem \ref{thm0}}

\noindent {\bf Proof of Theorem \ref{thm0} (1).} This  is a standard argument.
Since $u \in{{\Bbb X}_p}$, we have 
\begin{equation}
\int_\Omega |\nabla u|^{p-2}\nabla{u}\cdot \nabla\psi = \int_{\partial\Omega} \psi \nu+ \int_{\partial\Omega} \psi \mu \quad (\forall \psi \in C^1(\bar\Omega)),\label{adm}
\end{equation}
where $\mu=- \Delta_p u \in  \mathcal{M}_b(\Omega)$ and 
$\nu =|\nabla u|^{p-2}\frac{\partial u}{\partial n} \in \mathcal{M}_b(\partial\Omega)$. 
From  the  definition  we  have 
\begin{equation*}
[u]_{{{\Bbb X}_p}} =\sup_{\substack{\psi \in C^{1}(\bar{\Omega}) \\ \|{\psi}\|_{L^{\infty}} \le 1}} 
  \int_{\Omega}|\nabla u |^{p-2}\nabla{u}\cdot \nabla{\psi}\le \int_{\Omega} |{\Delta_p u}| + \int_{\partial\Omega}|\nabla u|^{p-2} \Bigl{|}{\frac{\partial{u}}{\partial{n}}}\Bigr{|}. \label{seminorm}
\end{equation*}
To see the  opposite inequality, without  the loss of generality  we assume  that  $\mu \in C^\infty_c(\Omega) $ and 
$\nu  \in C_c^{\infty}(\mathbb R^N)$  with $\supp \mu \cap \supp \nu =\phi$.
Define  $ \psi = \sgn (\mu) + \sgn({\nu})$, where $\sgn(t) = 1,  t>0; 0,  t=0 ; -1,  t<0$. Let $ \psi _\varepsilon$ be  a mollification of $\psi$ such  that
 $\psi_\varepsilon\in C^\infty_c(\mathbb R^N)$, $|\psi_\varepsilon|\le 1$ and $\psi_\varepsilon \to \psi$  as $\varepsilon\downarrow 0$.
Then for  any $\eta>0$ there exists a $\varepsilon>0$ such  that  we have 
$$    \int_{\Omega}|\nabla u |^{p-2}\nabla{u}\cdot \nabla{\psi_{\varepsilon}}\ge \int_{\Omega} |{\Delta_p u}| + \int_{\partial\Omega}|\nabla u|^{p-2} \Bigl{|}{\frac{\partial{u}}{\partial{n}}}\Bigr{|} -\eta.  $$
Since  $\eta  $ is  an arbitrary positive  number,  the desired inequality holds.\qed

\par\medskip
\par\noindent{\bf Proofs of  $(2)$ and $(3)$.}  The assertion (3)  clearly follows  from (2), we  hence  prove (2).
Assume  that   $u$  is admissible  in $\Bbb X_p$. Then  from Definition \ref{df1.2} there exists a sequence $\{ u_k\} \subset W^{1,p}(\Omega)\cap L^\infty(\Omega)$ satisfying the properties  $1$  and $2$. Noting  that $ \nabla (T_M u_k) = \chi_{|u_k|\le M}\nabla u_k$,  we have 
\begin{align*} 
  \int_\Omega\left| \nabla T_M(u_k)\right|^p \,dx&=  \int_\Omega\left| \nabla u_k\right|^{p-2}\nabla u_k\cdot \nabla T_M(u_k)\\
  &=\int_{\partial\Omega} |\nabla u_k|^{p-2} \frac{\partial u_k}{\partial n} T_M u_k- \int_\Omega  \Delta_p u_k T_M u_k\\
 & \le M [u_k]_{{{\Bbb X}_p}}.
\end{align*}
From the property  $1$ we  see that   $\Delta_p u_k \to \Delta_p u$ in $D'(\Omega)$  as  $k\to\infty$. From the  property $2$ together with  the  weak compactness of  Radon measures and  the uniqueness of weak  limit ( see also Remark \ref{remark2.1}.2 ),   
$\lim_{k\to\infty} \Delta_p u_k = \Delta_p u$ in the  sense of  measures.
Then by  Fatou's lemma the  assertion is    proved.\qed

\subsection{Proof of Theorem \ref{thm1}}

First we prove Theorem \ref{thm1} assuming that $u\in C^1(\overline\Omega)$ and $\Delta_p u\in L^1(\Omega)$.
Then  we treat the  general  case.

\begin{lem}\label{lem1}
Assume  that $u\in C^1(\bar\Omega)$ and $\Delta_p  u \in L^1(\Omega)$ (in the sense of distribution). Then
\begin{equation} 
\int_{\Omega} |\nabla u|^{p-2}\nabla{u^{+}}\cdot \nabla{\phi} \le \int_{\substack{\partial{\Omega} \\ [u \ge 0]}} \phi  |\nabla u|^{p-2} \frac{\partial{u}}{\partial{n}} - \int_{\substack{\Omega \\ [u \ge 0]}} \phi \Delta_p{u} \quad (\forall{\phi} \in C^{1}(\bar{\Omega}), \phi \ge 0 ~ in ~ \bar{\Omega})\label{4.13}
\end{equation}

\end{lem}
\noindent{\bf Proof.}
Let  $\Phi$ is a $C^2$ convex function in {$\mathbb R$}, $\Phi'\ge 0$ in $\mathbb R$ and $\Phi' \in L^\infty(\mathbb R)$. 
\par First we assume  that  $p\ge 2$. 

By a direct calculation we see that
 \begin{align} 
\Delta_p  \Phi(u)&  =   \Phi'(u)^{p-1} \Delta_p  u + 
  (p-1) \Phi'(u)^{p-2} \Phi''(u)| \nabla u|^p  \qquad  \mbox{ in } D'( \Omega).\label{5.2'} 
\end{align}
Since $ \Phi''\ge 0$,  we have
 \begin{align}\label{5d}
\Delta_p  \Phi(u)& \ge   \Phi'(u)^{p-1} \Delta_p  u  \qquad  \mbox{ in }  D'(\Omega),\end{align}
in particular,
$\Delta_p \Phi(u) \in L^1(\Omega)$.
Hence,  for any ${\phi} \in C^{1}(\bar{\Omega}), \phi \ge 0 ~ in ~ \bar{\Omega}$ we have
\begin{align} \label{5x}
\int_{\Omega} |\nabla \Phi(u)|^{p-2}\nabla{\Phi(u)}\cdot \nabla{\phi} 
& = \int_{\partial\Omega} |\nabla \Phi(u)|^{p-2}{\Phi'(u)}\frac{\partial u}{\partial n} {\phi} -\int_\Omega \Delta_p \Phi(u)\phi\\
&
\le \int_{\partial{\Omega} } \phi |\Phi'(u)|^{p-2} \Phi'(u)|\nabla u|^{p-2} \frac{\partial{u}}{\partial{n}} - \int_{\Omega} \phi |\Phi'(u)|^{p-2} \Phi'(u) \Delta_p{u}.\notag \end{align}
By the approximation  argument, this is  still valid for $C^1$ convex function $\Phi$. 
Now  we take a
 $\Phi$  in {$\mathbb R$} such that  $ \Phi(t)=t$  if $t \ge 0$, $|\Phi(t)|<1$ if $t<0$, $0\le \Phi' \le 1$ in {$ \mathbb R$} and $lim_{t\rightarrow -\infty} \Phi'(t)=0$. 
Set $\Phi_n(t)= \Phi(nt)/n$ for $t\in {\mathbb R}$ and $n=1,2,....$ Then we see that $\{\Phi_n\}$ is a sequence of $C^1$convex functions in
$ {\mathbb R} $ such that $ \Phi_n(t)=t$  if $t \ge 0$, $|\Phi_n (t)|<\frac{1}{n}$ if $t<0$, $0\le \Phi'_n \le 1$ in  { $\mathbb R$}. 
Then we see that  $ \Phi_n(t) \to t^+$ as $n\to\infty$.
Replacing  $\Phi$ by $\Phi_n$ in   (\ref{5x}) and letting $n\to\infty$,
 we  have  the desired inequality by the  dominated convergence  theorem.
\par We proceed to  the case where  $1<p<2$.  We set $\Phi^\eta(t):= \Phi(t) + \eta t$ for $t\in \mathbb R$ with $\eta >0$. Then we see  that for each $\eta >0$
\begin{equation}
\sup_{t\in \mbox{\bf \scriptsize$R$}} (\Phi^\eta)'(t)^{p-2}(\Phi^\eta)''(t)=\sup_{t\in{ \mbox{\bf \scriptsize$R$}}} (\Phi'(t)+\eta)^{p-2}\Phi''(t)\le \eta^{p-2}
\sup_{t\in \mbox{\bf \scriptsize$R$}} \Phi''(t)<\infty.
\end{equation}
Hence we can apply he previous argument with $\Phi^\eta$ instead of $\Phi$, so that in a similar way we reach to the inequality (\ref{5x}) replaced $\Phi$ by $\Phi^\eta$.
 Letting $\eta \rightarrow 0$, we have (\ref{4.13}) and 
this completes the proof.
\hfill$\Box$

\begin{lem}\label{lem2}
Assume  that $u\in C^1(\bar\Omega)$ and $\Delta_p  u \in L^1(\Omega)$ (in the sense of distribution). Then
$u^+ \in{{\Bbb X}_p}$ and 
\begin{equation} 
[u^+]_{{\Bbb X}_p} \le [u]_{{\Bbb X}_p} ~. \label{4.18}
\end{equation}
\end{lem}
\noindent{\bf  Proof.}
We  note  that $u^+ \in W^{1,p^*}(\Omega)$. For  the proof  of  Lemma it suffices to show the following.\begin{equation} 
\Big|{\int_\Omega|\nabla u|^{p-2} \nabla{u^+}\cdot \nabla\psi}\Big| \le [u]_{{\Bbb X}_p}\|\psi\|_{L^\infty} \quad (\forall \psi \in C^1(\bar\Omega)) ~.
\end{equation}
For $\tilde\psi \in C^1(\bar\Omega)$, we apply   (\ref{4.13}) with $\psi= \|\tilde\psi\|_{L^\infty}+\tilde\psi$. 
Then
\begin{equation} 
\begin{aligned}
\int_\Omega|\nabla u|^{p-2} \nabla{u^+}\cdot \nabla{\tilde\psi} \le \Big({\int_{\substack{\partial{\Omega} \\ [u \ge 0]}}|\nabla u|^{p-2}  \frac{\partial{u}}{\partial{n}} - \int_{\substack{\Omega \\ [u \ge 0]}} \Delta_p{u}}\Big)\|\tilde\psi\|_{L^\infty} \\ \qquad\qquad\qquad\qquad\qquad + \int_{\substack{\partial{\Omega} \\ [u \ge 0]}} \tilde\psi |\nabla u|^{p-2} \frac{\partial{u}}{\partial{n}} - \int_{\substack{\Omega \\ [u \ge 0]}} \tilde\psi \Delta_p{u}
\end{aligned}
\end{equation}
Noting that
\begin{equation*}
\int_{\substack{\partial{\Omega} \\ [u \ge 0]}} |\nabla u|^{p-2} \frac{\partial{u}}{\partial{n}} - \int_{\substack{\Omega \\ [u \ge 0]}} \Delta_p{u} = -\int_{\substack{\partial{\Omega} \\ [u < 0]}}|\nabla u|^{p-2}  \frac{\partial{u}}{\partial{n}} + \int_{\substack{\Omega \\ [u < 0]}} \Delta_p{u}
\end{equation*}
we have 
\begin{equation*}
\begin{aligned}
\int_\Omega |\nabla u|^{p-2} \nabla{u^+}\cdot \nabla{\tilde\psi} &\le - \Big({\int_{\substack{\partial{\Omega} \\ [u < 0]}} |\nabla u|^{p-2} \frac{\partial{u}}{\partial{n}} - \int_{\substack{\Omega \\ [u < 0]}} \Delta_p{u}}\Big)\|\tilde\psi\|_{L^\infty} + \int_{\substack{\partial{\Omega} \\ [u \ge 0]}} \tilde\psi |\nabla u|^{p-2}  \frac{\partial{u}}{\partial{n}} - \int_{\substack{\Omega \\ [u \ge 0]}} \tilde\psi \Delta_p{u} \\
&\le \Big({\int_{\partial{\Omega}} |\nabla u|^{p-2} \Big|\frac{\partial{u}}{\partial{n}}\Big| + \int_\Omega |\Delta_p{u}|}\Big)\|\tilde\psi\|_{L^\infty} = [u]_{{\Bbb X}_p}\|\tilde\psi\|_{L^\infty} ~.
\end{aligned}
\end{equation*}
By replacing $\tilde\psi$ by $-\tilde\psi$, we have  the desired inequality (\ref{4.18}).
\qed
\par\medskip
Secondly  we  assume  that $u$ is  admissible  in $\Bbb X_p$.
We recall a lemma on Neumann boundary problem for a monotone operator $\Delta_p$.

\begin{lem}\label{lem3} Let $\mu\in C^\infty_c(\Omega) $ and $\nu\in C_c^{\infty}(\Bbb R^N)$. Assume  that $\int _\Omega \mu + \int_{\partial\Omega}\nu =0$.\par
Then there exists a unique function  $u\in C^{1,\sigma} (\bar\Omega) $ for  some $\sigma\in (0,1)$ such that  
\begin{equation}\begin{cases} &-\Delta_p u =\mu\qquad \mbox { in  }\Omega \\  &|\nabla u|^{p-2} \frac{\partial u}{\partial n} =\nu \quad \mbox {on } \partial \Omega, \\
& \int_\Omega u  =0.
\end{cases}
\end{equation}

\end{lem}
\noindent{\bf  Proof.}  It follows from the standard theory 
that   we have the unique   solution  $u$ in  $W^{1,p}(\Omega)$.  For  the detail, refer to \cite{LL}; theorems 2.1 and  2.2  for example.
Since $ \mu$   and $\nu$ smooth, we see that 
 $u\in C^{1,\sigma} (\bar\Omega) $  for  some $\sigma \in (0,1)$ (See e.g.  DiBenedetto \cite{Di}).   Here we note  that $u$ is $p$-harmonic near the boundary as well. \qed
\vspace{1em}
\par
By Definition \ref{df1.2} of the admissibility in  $\mathbb X_p$ we  have for  each $k\ge 1$ that
\begin{equation}
\int_\Omega |\nabla u_k|^{p-2}\nabla{u_k}\cdot \nabla\psi = \int_{\partial\Omega} \psi|\nabla u_k|^{p-2} \frac{\partial{u_k}}{\partial{n}} - \int_{\partial\Omega} \psi \Delta_p{u_k} \quad (\forall \psi \in C^1(\bar\Omega)).\label{adm2}
\end{equation}
It  follows  from Remark \ref{remark2.1}(2) that  in the  sense of weak{*}  topology as $n\to\infty$
\begin{align}
\Delta_p{u_k} \stackrel{*}{\rightharpoonup} \Delta_p{u}  ~ in ~ \mathcal{M}_b({\Omega}),\quad & \|\Delta_p{u_k}\|_{L^1(\Omega)} \to \|\Delta_p{u}\|_{\mathcal{M}_b(\Omega)} ~,  \label{APP3}\\
|\nabla u_k|^{p-2} \frac{\partial{u_k}}{\partial{n}}\stackrel{*}{\rightharpoonup}|\nabla u|^{p-2} \frac{\partial{u}}{\partial{n}} ~ in ~ \mathcal{M}_b(\partial{\Omega}), \quad &  \||\nabla u_k|^{p-2}\frac{\partial{u_k}}{\partial{n}}\|_{L^1(\partial\Omega)} \to \Big{\|}|\nabla u|^{p-2}\frac{\partial{u}}{\partial{n}}\Big{\|}_{\mathcal{M}_b(\partial\Omega)} ~. \label{APP4}
\end{align}

By choosing $\psi=1$ in (\ref{adm2}), we have 
\begin{equation} \int_{\Omega}\Delta_p u_k = \int_{\partial\Omega} |\nabla u_k| ^{p-2}\frac{\partial u_k}{\partial n}. \label{4.23}\end{equation}
Let  us  set $\mu_k =-\Delta_p u_k$  and $\nu_k =|\nabla u_k|^{p-2} \frac{\partial u_k}{\partial n}$.
Let $\{\mu ^j_k\} \subset C^\infty_c (\bar\Omega)$ and $\{\nu^ j_k\}\subset C_c^\infty(\Bbb R^N)$ be  two sequences such that 
as $j\to\infty$
\begin{align}
\mu^{j}_k \stackrel{*}{\rightharpoonup} -\Delta_p{u_k} ~ weak^{*} ~ in ~  L^1({\Omega}),\quad & \quad \|\mu^j_k\|_{L^1(\Omega)} \to \|\Delta_p{u_k}\|_{L^1(\Omega)} ~,  \label{APP1}\\
\nu^{j}_k \stackrel{*}{\rightharpoonup}|\nabla u_k|^{p-2} \frac{\partial{u_k}}{\partial{n}} ~ weak^{*} ~ in ~ L^1(\partial{\Omega}), \quad & \quad \|\nu^j_k\|_{L^1(\partial\Omega)} \to \Big{\|}|\nabla u_k|^{p-2}\frac{\partial{u_k}}{\partial{n}}\Big{\|}_{L^1(\partial\Omega)} ~. \label{APP2}
\end{align}
From (\ref{4.23}) we  assume that
\begin{equation*}
\int_{\partial\Omega} \nu_k^j = -\int_\Omega \mu_k^j \quad ({}^\forall j,k \ge 1).
\end{equation*} 
It follows from Lemma \ref{lem3} that
for any $n\ge 1 $ and $ k \ge 1$, there exists  $w^k_n \in C^{1,\sigma}(\bar\Omega)$   such that
\begin{equation}
\begin{cases}
-\Delta_p{w_k^j} &= \mu_k^j \quad in ~ \Omega \\
|\nabla w_k^j|^{p-2}\frac{\partial{w_k^j}}{\partial{n}} &= \nu_k^j \quad on ~ \partial\Omega,
\end{cases}\label{app}
\end{equation}
or equivalently 
\begin{align} \int_{\Omega} |\nabla w^j_k|^{p-2}\nabla w^j_k \cdot\nabla \psi= \int_\Omega \psi\, d\mu^j_k +\int_{\partial\Omega}
\psi\, d\nu^j_k, \quad \mbox{ for  any } \psi\in C^1(\bar\Omega), \label{4.27}
\end{align}
 and without the loss of  generality  we also  assume that  for any $j,k \ge 1$
\begin{equation}
\int_\Omega w_k^j = \int_\Omega u_k. \label{average}
\end{equation}

Under these preparations we have 
\begin{lem}\label{lem4} For  each $n\ge 1$,   there exists a function $w_k\in W^{1,q}(\Omega)$ for every $q \in [1, \frac{N(p-1)}{N-1}) $ such that $w_k^j$  converges to $w_k$ in $w_k\in W^{1,q}(\Omega)$ as  $k\to\infty$ and $w_k$ satisfies (\ref{adm2}).
\end{lem}
\noindent{\bf  Proof.}
Since for each $k\ge 1$,  $\{\mu_k^j\}_{j=1}^\infty$  and $\{\nu_k^j\}_{j=1}^\infty$ are bounded in $ L^1(\Omega)$  and $ L^1(\partial\Omega)$ respectively,
this assertion follows  from the same argument in the proof of Theorem 1 in \cite{BG} with an obvious modification.
In fact,  one can show  that $\{ w^j_k\}_{j=1}^\infty$  is bounded in $W^{1,q}(\Omega)$, using  similar test functions  for $\psi$.
Then by the weak  compactness, Poincar\'e's inequality  and the Rellich type theorem, one can see that there exists a function $w_k\in W^{1,q}(\Omega)$ such that
 \begin{align*} &\nabla w^j_k \to \nabla w_k  \quad \mbox{ in } L^q \quad (\mbox{weak} )\\
&w^j_k \to  w_k  \quad \mbox{ in } L^q \\
 &w^j_k \to  w_k  \quad \mbox{ a.e. }  \end{align*}
Moreover one  can see that $\nabla w_k^j \to \nabla w_k$  in $L^1(\Omega)$.  Then by the dominated convergence theorem the conclusion follows in a quite similar way. For  the detail see \cite{BG}
\qed
\begin{lem}\label{lem6}  We have  $w_k= u_k$  a.e. for $k=1,2,\cdots$.

\end{lem}
\noindent{\bf  Proof.}
We claim that $w_k =  u_k\in W^{1,q}(\Omega)$ for  $q \in [1, \frac{N(p-1)}{N-1}) $.  Choose any  $M>0$.
Recalling that $u_k\in W^{1,p}(\Omega)\cap L^\infty(\Omega)$, we use $T_M(w^j_k-u_k )\in W^{1,p}(\Omega)\cap L^\infty(\Omega)$  as  a test function  in (\ref{adm2}) and (\ref{4.27}). By a  subtraction 
\begin{align}
 \int_\Omega (|\nabla w^j_k|^{p-2} \nabla w^j_k- & |\nabla u_k|^{p-2} \nabla u_k)  \cdot \nabla(T_M(w^j_k-u_k) \notag \\ &
=
 \int_{\Omega} T_M(w^j_k-u_k)\,d(\mu^j_k -\mu_k) +
 \int_{\partial\Omega} T_M(w^j_k-u_k)\,d(\nu^j_k -\nu_k). \notag
\end{align}
The left hand side is estimated from below in the following way,
\begin{equation}
 \int_\Omega (|\nabla w^j_k|^{p-2} \nabla w^j_k-  |\nabla u_k|^{p-2} \nabla u_k)  \cdot \nabla T_M(w^j_k-u_k)
\ge C\int_\Omega | \nabla T_M( w^j_k- u_k)|^p
\end{equation}
for some positive  number $C$ independent of  each $j$,
and  the right hand side  goes to $0$  as $j\to \infty$. 
Since  this  holds for  all $M>0$, we conclude by  the monotonicity of $\Delta_p$ that $\nabla w_k=\nabla u_k$ a.e. Taking into account that
$w_k \in  W^{1,q}(\Omega)$, $ u_k \in  W^{1,p}(\Omega)$ and (\ref{average}),
we conclude that $u_k=w_k \,$ a.e.
\qed
\par\medskip
\par\noindent{\bf End of proof of  Theorem \ref{thm1}.}
By applying Lemma \ref{lem2} we have 
\begin{equation} 
\Big|{\int_\Omega|\nabla (w_k^j)^+|^{p-2} \nabla{(w_k^j)^+}\cdot \nabla\psi}\Big| \le [w_k^j]_{{\Bbb X}_p}\|\psi\|_{L^\infty} \quad (\forall \psi \in C^1(\bar\Omega)) ~.
\end{equation}
From Lemma \ref{lem4} and Lemma \ref{lem6} we have, up to subsequence,  that $w_k^j \to u_k $ a.e. and $(w_k^j)_+ \to (u_k)_+$ in 
$W^{1,q}(\Omega)$ as  $j\to \infty$.
Letting $j \to \infty$, we have 
\begin{equation*}
\Big|{\int_\Omega|\nabla{u_k^+}|^{p-2} \nabla{u_k^+}\cdot \nabla\psi}\Big| \le [u_k]_{{\Bbb X}_p}\|\psi\|_{L^\infty} \quad (\forall \psi \in C^1(\bar\Omega)).
\end{equation*}
Finally letting $k\to\infty$ we have  the conclusion.
\qed

\subsection{Proof of Theorem \ref{thm2}}
\par \noindent {\bf Proof  of  the  assertion 1.}
\par\noindent{\bf 1st step.}
Assume  that $u$ is admissible in
$W^{1,p^*}_{0}(\Omega)$, and   hence  both $u^+$ and $u^-$ are admissible $W^{1,p^*}_{0}(\Omega)$.
From  the statement $4$ of Proposition \ref{prop1.1}, we   can assume  that 
 $\{ u_k\}_{k=1}^\infty \subset W^{1,p}_{0}(\Omega) \cap C^{1}_0(\Omega)$ in  Definition \ref{df1.1}.
We decompose $u_k \in W^{1,p}_{0}(\Omega) \cap C^{1}_0(\Omega)$   to  obtain 
$ u_k= u_k^+-u_k^-,\,\mbox{ where }\, u_k^+= \max(u_k, 0), \,  u_k^-= \max(-u_k,0). $
Then  each $ u_k^{\pm} \in  W^{1,p}_{0}(\Omega) \cap C_0^{1,0}(\bar\Omega)$. 
Since $u_k^+\ge 0$ in $\Omega$ and $u_k^+=0$ on $\partial \Omega$, we  see that $\frac{\partial u_k^+}{\partial n} \le 0$ on $\partial\Omega$.
Similarly we  have $\frac{\partial u_k^-}{\partial n} \le 0$ on $\partial\Omega$.
Therefore 
\begin{align*}- \int_{\partial\Omega} |\nabla u_k^+|^{p-2}\left|\frac{ \partial u_k^+}{\partial n}\right| &= 
 \int_{\partial\Omega} |\nabla u_k^+|^{p-2}\frac{ \partial u_k^+}{\partial n}=
 \int_\Omega \Delta_p u_k^+,\\
 - \int_{\partial\Omega} |\nabla u_k^-|^{p-2}\left|\frac{ \partial u^-}{\partial n}\right| &= 
 \int_{\partial\Omega} |\nabla u_k^-|^{p-2}\frac{ \partial u_k^-}{\partial n}=
 \int_\Omega \Delta_p u_k^-.
\end{align*}
Hence 
\begin{align*}
\int_{\partial\Omega} |\nabla u_k^+|^{p-2}\left|\frac{ \partial u_k^+}{\partial n}\right| \le 
\left| \int_\Omega \Delta_p u_k^+\right|,\quad
\int_{\partial\Omega} |\nabla u_k^-|^{p-2}\left|\frac{ \partial u_k^-}{\partial n}\right| \le 
\left| \int_\Omega \Delta_p u_k^-\right|.
\end{align*}
After all we have 
\begin{equation}
\int_{\partial\Omega} |\nabla u_k|^{p-2}\left|\frac{ \partial u_k}{\partial n}\right| \le 
 \int_\Omega |\Delta_p u_k|, \label{4.31}
\end{equation}
 in particular $ |\nabla u_k|^{p-2}\frac{ \partial u_k}{\partial n} \in L^1(\partial\Omega)$.
 Hence  we have
 \begin{equation} [u_k]_{{\Bbb X}_p}\le \int_{\partial \Omega} |\nabla u_k|^{p-2}\left|\frac{ \partial u_k}{\partial n}\right|+ \int_\Omega |\Delta_p u_k|\le 2 \int_\Omega |\Delta_p u_k|<\infty.
\end{equation}

\par\noindent{\bf 2nd step.} Again  assume  that $\{ u_k\}_{n=1}^\infty \subset W^{1,p}_{0}(\Omega) \cap C^{1}_0(\Omega)$ in  Definition \ref{df1.1}.
By Definition \ref{df1.1} $(1)$ we have 
\begin{equation} \int_\Omega |\nabla u_k|^{p-2} \nabla u_k\cdot \nabla \psi \to \int_\Omega |\nabla u|^{p-2} \nabla u\cdot \nabla \psi \quad \mbox{ for  any }  \psi\in C^1_c(\Omega).
\end{equation}
It  follows from the weak compactness of bounded measures and  the uniqueness of  weak limit  that 
 $\Delta_p u_k \to \Delta _p u$ strongly in $\mathcal M(\Omega)$. 
By  the previous step we have 
\begin{equation}|u_k|_{{\Bbb X}_p}\le 2 \int_\Omega |\Delta_p u_k|\qquad \mbox{  for  } k=1,2,\cdots.
\end{equation}
Hence  we see  that $ |\nabla u_k|^{p-2}\frac{\partial u_k}{\partial n} \in L^1(\partial \Omega)$ converge to some  measure $\nu$  in $ M(\partial\Omega)$ up to subsequences.
Therefore by  the lower semicontinuity of  the norm $||\cdot||_{M(\Omega)}$ with respect to the weak* convergence as $n\to\infty$, we have

$$ [u]_{{\Bbb X}_p}\le 2 \int_\Omega|\Delta_p u|.$$
Therefore $u$ is admissible in ${{\Bbb X}_p}$,  and  hence $u^+ \in {{\Bbb  X}_p}$ by Theorem \ref{thm1}. \hfill $\Box$
\par\bigskip
\par \noindent {\bf Proof  of  the  assertion 2.}
We claim  that $\int_{\Omega} |\Delta_p u^{+}| \le \int_\Omega |\Delta_p u|$.
\begin{lem}\label{lem7}
Assume  that $u \in C^1_0(\bar\Omega)$ and $\Delta_p u \in L^1(\Omega)$.  Then $\Delta{u^+} \in \mathcal{M}_b(\Omega)$ and 
\begin{equation} 
\|{\Delta{u^+}}\|_{\mathcal{M}_b(\Omega)} \le \|{\Delta{u}}\|_{L^1(\Omega)} ~.
\end{equation}
\end{lem}
\proof
By applying Lemma \ref{lem2} with $u+\varepsilon$, where $\varepsilon>0$, we deduce that
\begin{equation} |(u+\e)^+|_{{\Bbb X}_p}\le |u+\e|_{{\Bbb X}_p} = |u|_{{\Bbb X}_p}.
\end{equation}
Since $(u+\e)^+= u+\e $ in a nelghborhood of $\partial\Omega$, 
\begin{equation} \frac{\partial}{\partial n}(u+\e)^+= \frac{\partial u}{\partial n} \quad \mbox{ on } \partial\Omega.
\end{equation}
Noting that  
\begin{align*}|(u+\e)^+|_{{\Bbb X}_p} &= ||\Delta_p(u+\e)^+||_{\mathcal M(\Omega)}+ || |\nabla(u+\e)^+|^{p-2} \frac{\partial}{\partial n}(u+\e)^+||_{L^1(\partial\Omega)} \\
|u|_{{\Bbb X}_p} &= ||\Delta_p u||_{L^1(\Omega)} + || |\nabla u|^{p-2}  \frac{\partial u}{\partial n}||_{L^1(\partial\Omega)},
\end{align*}
 we immediately have 
\begin{equation}||\Delta_p(u+\e)^+||_{\mathcal M(\Omega)} \le || \Delta_p u||_{L^1(\Omega)} \quad \mbox{ for  any  } \e>0.
\end{equation}
The results follows from the lower semicontinuity  of  the  norm $ ||\cdot ||_{\mathcal M(\Omega)} $ with respect to the weak* convergence as $\e\to 0$.\qed

\subsection{Proof of  Theorem \ref{thm3}}
We prepare some fundamental lemmas.
\begin{lem}\label{lem8} Let  $u\in W^{1,p^*}(\Omega)$. Assume  that for some $h\in L^1(\partial \Omega)$  and $ g\in L^1(\Omega)$ we have 
\begin{equation}
\int_{\Omega} |\nabla u|^{p-2}\nabla u\cdot \nabla \varphi \le \int _{\partial \Omega } h \varphi  + \int _\Omega g \varphi  \quad 
\mbox{ for  any } \varphi \in C^1(\bar \Omega), \varphi\ge 0. \label{4.39}
\end{equation}
Then $u\in {{\Bbb X}_p}$. Moreover $-\Delta_p u \le g$  in $\mathcal M(\Omega) $  and $ |\nabla u|^{p-2}\frac{\partial  u}{\partial n}\le h$ in $\mathcal M(\partial\Omega)$.
\end{lem}
\noindent{\bf  Proof.}
 By (\ref{4.39}) we have 
 
\begin{equation}
\int_{\Omega} |\nabla u|^{p-2}\nabla u\cdot \nabla \varphi \le \int _{\partial \Omega } h^+ \varphi  + \int _\Omega g^+ \varphi  \quad 
\mbox{ for  any } \varphi \in C^1(\bar \Omega), \varphi\ge 0. \label{4.40}
\end{equation}
Using  nonnegative  test functions $||\varphi||_{L^\infty} \pm \varphi $ as  the  argument in the proof of Lemma \ref{lem2},   it  is  easy to see that
\begin{equation}
\left|  \int_{\Omega} |\nabla u|^{p-2}\nabla u\cdot \nabla \varphi \right| \le 2 ( ||h^+||_{L^1(\partial\Omega)}+|| g^+||_{L^1(\Omega)}) ||\varphi||_{L^\infty(\Omega)}.
\end{equation}
Then we see $u\in {{\Bbb X}_p}$. The rest of  the assertions are  clear.\qed

\begin{lem}\label{lem9} In the previous Lemma \ref{lem8}, we further assume  that $u$ is  admissible  in ${{\Bbb X}_p}$.  Then
we have 
\begin{equation} 
\int_{\Omega} |\nabla u^+|^{p-2} \nabla{u^{+}}\cdot \nabla{\varphi} \le \int_{\substack{\partial{\Omega} \\ [u \ge 0]}} h \varphi + \int_{\substack{\Omega \\ [u \ge 0]}} g \varphi \quad\mbox{ for  any } \varphi \in C^1(\bar \Omega), \varphi\ge 0.
\end{equation}
\end{lem}
By  the admissibility there exists a sequence $ \{ u_k\} \subset W^{1,p^*}(\Omega) $ having the  properties in Definition \ref{df1.2}. By virtue of Proposition \ref{prop1.1} we can assume  that $u_k \in C^1(\bar\Omega)$. Then it follows from Lemma \ref{lem1} that 

\begin{equation} 
\int_{\Omega} |\nabla u_k|^{p-2}\nabla{u_k^{+}}\cdot \nabla{\varphi} \le \int_{\substack{\partial{\Omega} \\ [u_k \ge 0]}} \varphi  |\nabla u_k|^{p-2} \frac{\partial{u_k}}{\partial{n}} - \int_{\substack{\Omega \\ [u_k \ge 0]}} \varphi \Delta_p{u_k} \quad (\forall{\varphi} \in C^{1}(\bar{\Omega}), \varphi \ge 0 ~ in ~ \bar{\Omega})
\end{equation}
Taking a limit as $k\to \infty$ we have  
\begin{equation} 
\int_{\Omega} |\nabla u|^{p-2}\nabla{u^{+}}\cdot \nabla{\varphi} \le \int_{\substack{\partial{\Omega} \\ [u \ge 0]}} \varphi  |\nabla u|^{p-2} \frac{\partial{u}}{\partial{n}} - \int_{\substack{\Omega \\ [u\ge 0]}} \varphi \Delta_p{u} \quad (\forall{\varphi} \in C^{1}(\bar{\Omega}), \varphi \ge 0 ~ in ~ \bar{\Omega})
\end{equation}
Using  Lemma \ref{lem6} the conclusion holds.\qed

\begin{lem}\label{lem10}

Assume  that $u \in C^1(\bar\Omega)$ is  admissible in ${{\Bbb X}_p}$ and 

$|\nabla u|^{p-2}\frac{\partial{u}}{\partial{n}} \in L^{1}(\partial{\Omega})$. Then
\begin{align} 
|\nabla u^+|^{p-2}\frac{\partial{u^{+}}}{\partial{n}} \le \begin{cases}
							|\nabla u|^{p-2}\frac{\partial{u}}{\partial{n}} \quad &\mbox{on }~ [u > 0] \\
							0 \quad &\mbox{on } ~ [u < 0] \\
							\min \{ |\nabla u|^{p-2}\frac{\partial{u}}{\partial{n}}, 0 \} &\mbox{on } ~ [u = 0].
						    \end{cases}
\end{align}

\end{lem}
\noindent{\bf  Proof.}
Put
$\mu = (-\Delta_p{u})^{+},~ h = |\nabla u|^{p-2}\frac{\partial{u}}{\partial{n}}$. Then
\begin{equation*}
\int_{\Omega} |\nabla u|^{p-2} \nabla u \cdot \nabla \varphi  \le \int_{\partial{\Omega}} h \varphi + \int_{\Omega} \varphi~d\mu \quad (\forall{\varphi} \in C^{1}(\bar{\Omega}),\varphi \ge 0 ~ in ~ \bar{\Omega})
\end{equation*}

It follows from Lemma \ref{lem9} that  $u^{+}$ satisfies
\begin{equation} 
\int_{\Omega} |\nabla u|^{p-2}\nabla{u^{+}}\cdot \nabla{\varphi} \le \int_{\substack{\partial{\Omega} \\ [u \ge 0]}} h \varphi + \int_{\Omega} \varphi~d\mu \quad (\forall{\varphi} \in C^{1}(\bar{\Omega}),~ \varphi \ge 0 ~ in ~ \bar{\Omega})
\end{equation}
By Theorem \ref{thm1}  we  have $u^{+} \in{{\Bbb X}_p}$, hence
\begin{equation}
|\nabla u|^{p-2}\frac{\partial{u^{+}}}{\partial{n}} \le \chi_{[u \ge 0]} h  = \chi_{[u \ge 0]} |\nabla u|^{p-2}\frac{\partial{u}}{\partial{n}} \quad on ~ \partial{\Omega}.
\end{equation}
By using $u-\varepsilon$, where $\varepsilon>0$ instead of $u$ we have in a similar way  that

\begin{equation}
|\nabla u|^{p-2}\frac{\partial{u^{+}}}{\partial{n}} \le \chi_{[u >0]} h  = \chi_{[u > 0]} |\nabla u|^{p-2}\frac{\partial{u}}{\partial{n}} \quad on ~ \partial{\Omega}.
\end{equation}
In particular,
\begin{equation} 
|\nabla u|^{p-2}\frac{\partial{u^{+}}}{\partial{n}} \le 0 \quad on ~ [u = 0]
\end{equation}
Hence the conclusion follows.
\hfill$\Box$

\begin{cor}
Assume  that $u$ is  admissible in ${{\Bbb X}_p}$  and 
$u \in  W^{1,p^*}_{0}(\Omega)$. If $u \ge 0 ~ in ~ \Omega$, then 
\begin{equation*}
|\nabla u|^{p-2}\frac{\partial{u}}{\partial{n}} \le 0 \quad on ~ \partial{\Omega}.
\end{equation*}

\end{cor}
\noindent{\bf  Proof.}
$u = u^{+} ~ in ~ \Omega$  and A$u = 0 ~ on ~ \partial{\Omega}$, hence applying the Lemma \ref{lem10} we have A
\begin{equation*}
\frac{\partial{u}}{\partial{n}} = \frac{\partial{u^{+}}}{\partial{n}} \le \min \{ \frac{\partial{u}}{\partial{n}}, 0 \} \le 0 \quad on ~ \partial{\Omega}.
\end{equation*}
\hfill$\Box$

\noindent{\bf  Proof of Theorem \ref{thm3}.}
By Theorem \ref{thm1} $u^{+} \in{{\Bbb X}_p}$.  By applying Kato's inequality (   Corollary 1.1 in \cite{LHkato}  ) to  $u-a \in {{\Bbb X}_p}$,
we havre
\begin{equation*}
\Delta_p{(u-a)^{+}} \ge \chi_{[u \ge a]} \Delta_p{u} = G \quad in ~ \Omega
\end{equation*}
for  any $a\in  {\bf R}$. Here we  note  thatt $ (\Delta_p u)_d=\Delta_p u $, because  $\Delta_p u \in L^1(\Omega)$.
 Letting $a \downarrow 0$ we  have  
\begin{equation*}
\Delta_p{u^{+}} \ge \chi_{[u > 0]} \Delta_p{u} = G \quad in ~ \Omega.
\end{equation*}
Combining  this with Lemma \ref{lem8},  we have  for any $\varphi \in C^{1}(\bar{\Omega}),~ \varphi \ge 0 ~ in ~ \Omega$
\begin{equation*}
\int_{\Omega} |\nabla u|^{p-2}\nabla{u^{+}}\cdot \nabla{\varphi} = \int_{\partial{\Omega}} \varphi |\nabla u|^{p-2} \frac{\partial{u^{+}}}{\partial{n}} - \int_{\Omega} \varphi \Delta{u^{+}} \le \int_{\partial{\Omega}} H \varphi - \int_{\Omega} G \varphi.
\end{equation*}
\hfill$\Box$

\section{Appendix ( Proposition \ref{prop1.1})}
We begin with  recalling  a local version of Admissibility in \cite{LH}.
  \begin{df}\label{df3} {\bf  $($Admissibility   in $W^{1,p^*}_{\rm loc}(\Omega)$$)$}  Let  $ 1<p<\infty $ and $p^* =\max[1,p-1]$.
 A function $u$ is said to be  admissible in   in $W^{1,p^*}_{\rm loc}(\Omega)$, if   $u \in W^{1,p^*}_{\rm loc}(\Omega)$, $\Delta_p u \in \mathcal{M}(\Omega)$ ; the total measure  is not necessarily  finite,  and if     
 there exists a sequence $\{ u_k\}_{k=1}^\infty \subset W^{1,p}_{\rm loc}(\Omega) \cap L^\infty(\Omega)$  such  that: 
 \begin{enumerate}
    \item  $u_k\to u $ a.e. in $\Omega$ and $\,  u_k\to u $ in $ W^{1,p^*}_{\rm loc}(\Omega)$ as $k \to \infty$.
 \item $\Delta_p u_k\in L^1_{\rm loc}(\Omega)\,$  $(k=1,2,\cdots)$ and 
 \begin{equation} \sup_{k}|\Delta_p u_k|(\omega)= \sup_{k}\int_\omega |\Delta_p u_k|<\infty \quad \mbox{  for all } \omega\subset \subset \Omega.
 \end{equation}
  \end{enumerate}
  \end{df}
Here we describe the following fundamental results, parts of  which are  already known. 
 \begin{prop}  \label{prop1.1}
 \begin{enumerate}
 Let   
 $\Omega$  be a bounded smooth  domain of {\bf $\mathbb R^{ N}$}. 
\item  Assume  that    $u$ is   admissible in $W_{\rm loc}^{1,p^*}(\Omega)$.  Then, for  every $M>0$,
   $T_M u\in W_{\rm loc}^{1,p}(\Omega) $.
\item 
 A  function $u\in W^{1,p}_{0}(\Omega)$  is  admissible in $W^{1,p^*}_{0}(\Omega)$, if    $\Delta_p u \in \mathcal{M}_b(\Omega)$.
 \item 
 A  function $u\in W^{1,p}_{\rm loc}(\Omega)$  is  admissible in $W^{1,p^*}_{\rm loc}(\Omega)$, if    $\Delta_p u \in \mathcal{M}(\Omega)$.

 \item In Definition \ref{df1.2}, the sequence $\{u_k\} $ can be taken in $C^1(\bar\Omega)$. 
  \item In Definition \ref{df1.1}, the sequence $\{u_k\} $ can be taken in $C^{1}_0(\bar\Omega)=\{\varphi\in C^{1}(\bar\Omega): u=0 \mbox{ on } \partial\Omega\}$.  

 \end{enumerate}

   \end{prop}

 The proof  of assertion $1$  for $p=2$ is seen 
 in  \cite{BP1} and \cite{BP2}) and for   $p>1$  in \cite{LH},  and  the proof  of  assertion $2$  is   seen in Appendix of \cite{LH}. The assertion $4$ is already verified in  the proof  of Theorem \ref{thm1}.
Therefore we establish the  assertions $3$ and  $5$ in  the  rest of  this  section.

   \par\medskip\noindent{\bf Proof of  assertion 3.} 
To use a diagonal argument, we choose and fix  a family of open set $\{ \omega_k\}$ such that
\begin{equation} \omega_1\subset\subset  \omega_2\subset\subset \cdots \subset\subset  \omega_k \subset\subset  \omega_{k+1} \subset\subset\cdots \subset\subset\Omega
\mbox{  and } \Omega = \cup_{k=0}^\infty \omega_k.\label{diagonal} \end{equation} 
  Let $\rho \in C_{0}^\infty(B_1)$ be a radial,
 nonnegative and decreasing mollifier. By extending $ v\in L^1(\Omega)$ to the whole space  so that $v \equiv 0$  outside $\Omega$,  we define a mollification of $v$  with $\varepsilon>0$ by
 \begin{equation}
  v^\varepsilon(x):=\rho_{\varepsilon} \ast v(x) = \int_\Omega \rho_{\varepsilon} (x-y)v(y)dy  \qquad {\rm for}\,\, x \in \Omega. \label {mollifier}
 \end{equation}
First we  prove that
 $u\in W^{1,p}_{0}(\Omega)$  is  admissible in  $W^{1,p^*}_{\rm loc}(\Omega)$, if   $\Delta_p u$ is  a Radon measure on $\Omega$. 
 Again by extending $u \in W^{1,p}_{0}(\Omega) $ and $\Delta_p u\in W^{-1, p'}$ to the whole space so that  $ u=0$  and $\Delta_p u =0$ outside $\Omega$ respectively. 
 Let $w_k \in W^{1,p}_0(\Omega) \cap C^1(\overline\Omega)$  be the unique  weak solution  of the boundary  value problem  for  the monotone operator $\Delta_p$ (see e.g. \cite {LL}): For  $k=1,2,\cdots$ and $\varepsilon_1>\varepsilon_2>\cdots \varepsilon_k>\cdots \to 0$, we set
 \begin{equation}\label{4.3} \begin{cases}& \Delta_p w_k = (\Delta_p u)^{\varepsilon_k} \qquad \mbox{ in } \Omega,\\
 & w_k =0 \qquad \mbox{ on } \partial \Omega,
 \end{cases}
 \end{equation}
where  
$ |\nabla u|^{p-2}\nabla u\in (L^{p'}(\Omega))^N$ with $p'=p/(p-1)$, $(|\nabla u|^{p-2}\nabla u)^{\varepsilon_k}\in (C^\infty(\mathbb R^N))^N $
and  $(|\nabla u|^{p-2}\nabla u)^{\varepsilon_k}$ is a mollification of $|\nabla u|^{p-2}\nabla u$ defined by (\ref{mollifier}).
Let us set $\Delta_pu= \mu$.  We note  that  $|\mu|(\omega)<\infty$ for any $\omega\subset\subset \Omega$. 
Then we have    $ \div (|\nabla u|^{p-2}\nabla u)^{\varepsilon_k}= (\div |\nabla u|^{p-2}\nabla u)^{\varepsilon_k} = (\Delta_p u)^{\varepsilon_k}=\mu^{\varepsilon_k}$ in $\omega$ provided that  $\varepsilon_k$ is sufficiently small.
 Hence  we  clearly have
 \begin{align*}
&| \Delta_p w_k |(\omega)= |\mu^{\varepsilon_k}|(\omega) \to |\mu|(\omega) \,\, \mbox{ as } \,\,  k\to\infty.
 \end{align*}
Since $\mu$ does  not  charge $\partial\Omega$, this proves  the condition 2.  Next we show 
  \begin{equation}\label{7.2}
 w_k \to u \mbox{ in } W^{1,p}_{0}(\Omega)  \,\, \mbox{ as } \,\,  k\to\infty.\end{equation} 
Then we can choose a subsequence so  that the  condition 1 is   satisfied.
 By using $w_k-u \in W^{1,p}_{0}(\Omega) $  as  a test function, we have 
 \begin{align}\label{7.9}
 - \langle \Delta_p w_k- \Delta_p u, w_k-u \rangle &= \int_\Omega  |(\nabla w_k|^{p-2}\nabla w_k - |\nabla u|^{p-2}\nabla u)\cdot \nabla(w_k-u) \notag \\
 &\ge c_2\int_\Omega | \nabla( w_k- u)|^p.
  \end{align}
In  the left-hand side,   using  Young's inequality for $\delta>0$ we have 
 \begin{align}\label{7.10}
 -&\langle \Delta_p w_k-\Delta_p u, w_k-u \rangle 
= \int_\Omega( (|\nabla u|^{p-2}\nabla u )^{\varepsilon_k}- |\nabla u|^{p-2}\nabla u)\cdot \nabla(w_k-u) \notag \\
 & \le C(\delta)\int_\Omega  | (|\nabla u|^{p-2}\nabla u )^{\varepsilon_k}-  |\nabla u|^{p-2}\nabla u|^{p'}+ \delta \int_\Omega |\nabla(w_k-u)|^p ,
 \end{align}
 where $C(\delta)>0$ is a constant depending only on $\delta$.\par
 We note  that   $||\,(|\nabla u|^{p-2}\nabla u )^{\varepsilon_k}-  |\nabla u|^{p-2}\nabla u \, ||_{L^{p'}(\Omega)}\to 0$ as $k\to\infty$.
 It follows from (\ref{7.9}) and (\ref{7.10}) that $\nabla w_k \to \nabla u$  in $(L^{p}(\Omega))^N$ as $n\to \infty$, which implies (\ref{7.2}).
 Then, taking a subsequence if  necessary, $\{w_k\} \subset  W^{1,p}_0(\Omega) \cap C^1(\overline\Omega)$ 
 satisfies the property $w_k\to u$  a.e.  in $\Omega$ as $k\to \infty$.\par
 Lastly we treat the case where $u\in W^{1,p}_{\rm loc}(\Omega)$. For each $k$ we choose $\eta_k\in C_c^\infty(\omega_{k+1})$ such that  $0\le \eta_k \le 1$ and
$\eta_k =1 $ in some neighborhood of $\overline{\omega_{k}}$. Let us set $ v_k= \eta_k u\,(k=1,2,3,\cdots)$ . Then
we see that  $v_k \in W^{1,p}_0 (\omega_{k+1})$, $ v_k \to u$ in $W^{1,p}_{\rm loc}(\Omega) $ as $k\to\infty$ and  $\Delta_pv_k \in W^{-1, p'}(\Omega)\cap M_b(\omega_k)$. Moreover   we have 
$|\Delta_p v_k|({\omega_j})=|\Delta_p u|(\omega_j)$ for any  $k\ge j$. Hence 
 $u$ is admissible in $W^{1,p^*}_{\rm loc}(\omega_k)$ with  $\Delta_p u\in \mathcal M_b(\omega_k)$  having an admissible sequence $\{v_k\}$. 
 By  the previous step with obvious modification, one can approximate each $v_k$ inductively by $\xi_k \in W_0^{1,p}(\Omega)\cap C^1(\overline\Omega)$  such  that $\xi_k\to u$ in $W^{1,p^*}_{\rm loc}(\Omega)$  as $k\to\infty$ and 
$||\Delta_p \xi_k|(\omega_j) - |\Delta_p u|(\omega_j) |<\frac 1k$ for $k\ge j$. Therefore  the assertion is now  proved.
\hfill$\Box$

\par\medskip\noindent 
{\bf Proof of  assertion $5$.} We  assume  that $u$ is  admissible  in $W^{1,p^*}_{0} (\Omega)$. Then we have a sequence of functions $\{u_k\}\subset 
W^{1,p}_{0} (\Omega)\cap L^{\infty}(\Omega)~(k=1,2,\ldots ) $  satisfying the properties 1 and 2 in Definition \ref{df1.1}.
By  the previous  step,  we  see that each $u_k$ is approximated  as  $j\to\infty$ by a sequence of  functions $\{w_k^j\} \subset  W^{1,p}_0(\Omega) \cap C^1(\overline\Omega)$ defined by
(\ref{4.3}) with $w_k= w_k^j$, $ u=u_k$  and $\varepsilon_k = \varepsilon_j$. 
Then we choose 
 a suitable subsequence of $\{ w_k^{j_k}\}$ as an approximation of $u$  so that  the assertion is  verified.
\hfill$\Box$
%
%
%
%

Toshio Horiuchi:\par
Department of Mathematics, Faculty of Science, Ibaraki University \par
Bunkyo 2-1-1 Mito City Ibaraki Prefecture  Japan
e-mail;  thoshio.horiuchi.math@vc.ibaraki.ac.jp

\par\medskip

Peter Kumlin
\par
Department of Mathematics, Chalmers University of Technology \par
S-41296 Goteborg, Sweden \par
e-mail; kumlin@chalmers.se

\end{document}